\documentclass{article}

\usepackage{graphicx}

\newtheorem{thm}{Theorem}[section]
\newtheorem{cor}[thm]{Corollary}
\newtheorem{lem}[thm]{Lemma}
\newtheorem{defi}[thm]{Definition}
\newcounter{bean}
\newcounter{milk}

\usepackage{graphicx}
\begin{document}

\title{{\bf The Mystery of the Shape Parameter II}}         
\author{{\bf Lin-Tian Luh} \thanks{supported by the MOST project 107-2115-M-126-005- } \\Department of Mathematics, Providence University \\Shalu, Taichung, Taiwan\\email:ltluh@pu.edu.tw\\phone:(04)26328001 ext. 15126\\ fax:(04)26324653 }        
\date{\today}          
\maketitle

\bigskip
{\bf Abstract.} We continue an earlier study of the shape parameter c contained in the famous multiquadrics $(-1)^{\lceil \beta \rceil}(c^{2}+\|x\|^{2})^{\beta},\ \beta>0$, and the inverse multiquadrics $(c^{2}+\|x\|^{2})^{\beta},\ \beta<0$. Instead of the space of bandlimited functions, we are going to treat a more general function space which roughly speaking is the same as the native space of gaussians. A totally different set of criteria for the optimal choice of c will be provided.
\\
{\bf keywords}: radial basis function, multiquadric, shape parameter, interpolation.
\section{{\bf Introduction}}
\hspace{5mm}As before, we are going to adopt a seemingly more complicated definition
\begin{equation}
h(x):=\Gamma \left (-\frac{\beta}{2}\right )(c^{2}+\|x\|^{2})^{\frac{\beta}{2}},\ \beta \in R \backslash 2N_{\geq 0},\ c>0,
\end{equation}
where $\|x\|$ is the Euclidean norm of $x$ in $R^{n}$, $\Gamma$ is the classical gamma function, and $c,\ \beta$ are constants. This definition will relieve our pain of manipulating its Fourier transform and developing useful criteria.

Recall that $h(x)$ is conditionally positive definite (c.p.d.) of order $m=max\left \{ \left \lceil \frac{\beta}{2} \right \rceil ,\ 0\right \}$ where $\left \lceil \frac{\beta}{2}\right \rceil$ denotes the smallest integer greater than or equal to $\frac{\beta}{2}$. This will be used in the text.

For the reader's convenience we review some basic features of the development in \cite{Lu4}. For any interpolated function $f$, the interpolating function will be of the form
\begin{equation}
s(x):=\sum_{i=1}^{N}c_{i}h(x-x_{i})+p(x)
\end{equation}
where $p(x)\in P_{m-1}$, the space of polynomials of degree less than or equal to $m-1$ in $R^{n}$, $X=\{x_{1},\cdots ,x_{N}\}$ is the set of centers (interpolation points). For $m=0,\ P_{m-1}:=\{0\}$. We require that $s(\cdot )$ interpolate $f(\cdot )$ at data points $(x_{1},f(x_{1})), \cdots , (x_{N},f(x_{N}))$. This leads to a linear system of the form
\begin{eqnarray}
\sum_{i=1}^{N}c_{i}h(x_{j}-x_{i})+\sum_{i=1}^{Q}b_{i}p_{i}(x_{j})=f(x_{j})&  & ,\ j=1,\cdots , N \nonumber \\
                                                                          &  &  \\
\sum_{i=1}^{N}c_{i}p_{j}(x_{i})=0 \hspace{3.45cm} &                                          & ,\ j=1,\cdots , Q \nonumber 
\end{eqnarray}
to be solved, where $\{ p_{1},\cdots,p_{Q}\}$ is a basis of $P_{m-1}$.

The sovability of the linear system is guaranteed by the c.p.d. property of $h$. However if $c$ is very large, $h$ will be numerically constant, making the linear system numerically unsolvable. Moreover, as pointed out by Madych in \cite{MN3}, if $c$ is very large, the coefficient matrix of the linear system will have a very large condition number, making the interpolating function $s$ unreliable when $f(x_{1}),\cdots ,f(x_{N})$ are not accurately evaluated.

Each function of the form (1) induces a function space called  native space denoted by $ {\cal C}_{h,m}(R^{n})$, abbreviated as $ {\cal C}_{h,m}$, where $m$ denotes its order of conditional positive definiteness. For its definition and characterization we refer the reader to \cite{Lu1},\cite{Lu2},\cite{MN1},\cite{MN2} and \cite{We}. This space is closely related to our space of interpolated functions and has to be used.

As in \cite{Lu4}, we need the following basic definitions for our development of the criteria.
\begin{defi}
For $n=1,2,3,\cdots ,$ the sequence of integers $\gamma_{n}$ is defined by $\gamma_{1}=2$ and $\gamma_{n}=2n(1+\gamma_{n-1})$ if $n>1$.
\end{defi}
\begin{defi}
Let $n$ and $\beta$ be as in (1). The numbers $\rho$ and $\Delta_{0}$ are defined as follows.
\begin{list}
  {(\alph{bean})}{\usecounter{bean} \setlength{\rightmargin}{\leftmargin}}
  \item Suppose $\beta <n-3$. Let $s=\left \lceil \frac{n-\beta -3}{2}\right \rceil $. Then 
    \begin{list}{(\roman{milk})}{\usecounter{milk} \setlength{\rightmargin}{\leftmargin}}
      \item if $\beta <0,\ \rho=\frac{3+s}{3}\ and\  \Delta_{0}=\frac{(2+s)(1+s)\cdots 3}{
 \rho^{2}};$
      \item if $\beta >0,\ \rho=1+\frac{s}{2\left \lceil \frac{\beta}{2}\right \rceil +3} \ and \ \Delta_{0}=\frac{(2m+2+s)(2m+1+s)\cdots (2m+3)}{\rho^{2m+2}}$ \\
where $ m=\left \lceil \frac{\beta}{2}\right \rceil$.          
    \end{list}
  \item Suppose $n-3\leq \beta <n-1$. Then $\rho=1$ and $\Delta_{0}=1$.
  \item Suppose $\beta \geq n-1$. Let $s=-\left \lceil \frac{n-\beta -3}{2}\right \rceil $. Then
 $$\rho =1\ and \ \Delta_{0}=\frac{1}{(2m+2)(2m+1)\cdots (2m-s+3)} \ where \ m=\left \lceil \frac{\beta}{2}\right \rceil.$$  
\end{list}
\end{defi}

The following theorem is cited directly from \cite{Lu4}.
\begin{thm}
  Let $h$ be defined as in (1) and $m=max\left \{ 0, \left \lceil \frac{\beta}{2}\right \rceil \right \}$. Then given any positive number $b_{0}$, there are positive constants $\delta_{0}$ and $\lambda,\ 0<\lambda <1$, which depend completely on $b_{0}$ and $h$ for which the following is true: For any cube $E$ in $R^{n}$ of side length $b_{0}$, if $f\in {\cal C}_{h,m}$ and $s$ is the map defined as in (2) which interpolates $f$ on a finite subset $X$ of $E$, then
\begin{equation}
  |f(x)-s(x)| \leq 2^{\frac{n+\beta+1}{4}}\pi^{\frac{n+1}{4}}\sqrt{n\alpha_{n}}c^{\frac{\beta}{2}}\sqrt{\Delta_{0}}(\lambda)^{\frac{1}{\delta}}\| f\| _{h}
\end{equation} 
holds for all $0<\delta \leq \delta_{0}$ and all $x$ in $E$ provided that $\delta=d(E,X):=\sup _{y\in E}\inf _{x\in X}|y-x|$. Here, $\alpha_{n}$ denotes the volume of the unit ball in $R^{n}$, and $c,\ \Delta_{0}$ were defined in (1) and Definition 1.2 respectively. Moreover $\delta_{0}=\frac{1}{6C\gamma_{n}(m+1)}$, and $\lambda=(\frac{2}{3})^{\frac{1}{6C\gamma_{n}}}$ where
$$C=\max \left\{ 2\rho'\sqrt{n}e^{2n\gamma_{n}},\ \frac{2}{3b_{0}}\right\},\ \rho'=\frac{\rho}{c}.$$
The integer $\gamma_{n}$ was defined in Definition 1.1, and $\| f\| _{h}$ is the $h$-norm of $f$ in ${\cal C}_{h,m}$. The constant $\rho$ was defined in Definition 1.2.
\end{thm}
{\bf Remark}.
Obviously the domain $E$ in Theorem 1.3 can be extended to a more general set $\Omega \subseteq R^{n}$ which can be expressed as the union of rotations and translations of a fixed cube of side $b_{0}$.
\\

In this paper the space of interpolated functions is defined as follows.
\begin{defi}
For any positive number $\sigma$,
$$E_{\sigma}:=\{f\in L^{2}(R^{n}):\ \int |\hat{f}(\xi)|^{2}e^{|\xi|^{2}/\sigma}d\xi <\infty \}$$
where $\hat{f}$ denotes the Fourier transform of $f$.
\end{defi}    
{\bf Remark}. It's easily seen that $E_{\sigma}$ is just the well-known native space of gaussian. For each $f$ in $E_{\sigma}$, we define its norm by
$$\|f\|_{E_{\sigma}}:=\left\{ \int|\hat{f}(\xi)|^{2}e^{|\xi|^{2}/\sigma}d\xi \right\} ^{1/2}.$$
\section{{\bf Fundamental Theory}}

\hspace{5mm}It's easily seen from Theorem 1.3 that the error bound (4) is greatly influenced by the choice of $c$. This is indeed the starting point of our theory. However, in order to develop useful criteria for the choice of $c$, some technical manipulation and theoretical analysis are necessary.
\begin{lem}
Let $\sigma>0$ and $\beta<0$. If $|n+\beta|\geq 1$ and $n+\beta+1\geq 0$, then $E_{\sigma}\subseteq {\cal C}_{h,m}(R^{n})$ and for any $f\in E_{\sigma}$,
$$\|f\|_{h}\leq 2^{-n-\frac{1+\beta}{4}}\pi^{-n-\frac{1}{4}}c^{\frac{1-n-\beta}{4}}\left\{ \sup_{|\xi|\in R^{n}}|\xi|^{(n+\beta+1)/2}e^{c|\xi|-|\xi|^{2}/\sigma}\right\}^{1/2}\|f\|_{E_{\sigma}}$$
where $\|f\|_{h}$ is the $h$-norm of $f$ in the native space ${\cal C}_{h,m}(R^{n})$.
\end{lem}
{\bf Proof}. Let $f\in E_{\sigma}$. By \cite{MN2} and \cite{Lu1},
\begin{eqnarray*}
\|f\|_{h} & = & \left\{ \frac{1}{(2\pi)^{2n}}\int\frac{|\hat{f}(\xi)|^{2}}{\hat{h}(\xi)}d\xi\right\}^{1/2}\\
          & =  & \left\{\frac{1}{(2\pi)^{2n}}\int\frac{|\hat{f}(\xi)|^{2}}{2^{1+\frac{\beta}{2}}\left(\frac{|\xi|}{c}\right)^{-\frac{\beta}{2}-\frac{n}{2}}{\cal K}_{\frac{n+\beta}{2}}(c|\xi|)}d\xi\right\}^{1/2}\ (Theorem\ 8.15\ of\ \cite{We}) \\
          & \leq & \left\{\frac{1}{(2\pi)^{2n}}\int \frac{|\hat{f}(\xi)|^{2}\sqrt{c|\xi|}e^{c|\xi|}}{2^{1+\frac{\beta}{2}}\left(\frac{|\xi|}{c}\right)^{-\frac{\beta+n}{2}}\sqrt{\frac{\pi}{2}}}d\xi\right\}^{1/2}\ (Corollary\ 5.12\ of\ \cite{We})\\
          & = & \frac{c^{\frac{1-n-\beta}{4}}}{2^{n+\frac{1+\beta}{4}}\pi^{n+\frac{1}{4}}}\left\{\int|\hat{f}(\xi)|^{2}|\xi|^{\frac{n+\beta+1}{2}}e^{c|\xi|}d\xi\right\}^{1/2}\\
          & \leq & \frac{c^{\frac{1-n-\beta}{4}}}{2^{n+\frac{1+\beta}{4}}\pi^{n+\frac{1}{4}}}\left\{\sup_{|\xi|\in R^{n}}\frac{|\xi|^{\frac{n+\beta+1}{2}}e^{c|\xi|}}{e^{\frac{|\xi|^{2}}{\sigma}}}\right\}^{1/2}\|f\|_{E_{\sigma}}\\
          & < & \infty.
\end{eqnarray*}
Hence  $E_{\sigma}\subseteq {\cal C}_{h,m}(R^{n})$  by Corollary 3.3 of \cite{MN2} and the lemma follows. \hspace{4.4cm} $\sharp$
\\

In the preceding proof we didn't find the supremum of the function in the braces. Let's try it. Suppose
$$G(\xi):=\xi^{\frac{n+\beta+1}{2}}e^{c\xi-\frac{\xi^{2}}{\sigma}},\ \xi>0.$$ Then
\begin{eqnarray*}
G'(\xi) & = & \frac{n+\beta+1}{2}\xi^{\frac{n+\beta-1}{2}}e^{c\xi-\frac{\xi^{2}}{\sigma}}+\xi^{\frac{n+\beta+1}{2}}e^{c\xi-\frac{\xi^{2}}{\sigma}}\left(c-\frac{2}{\sigma}\xi\right)\\
        & = & e^{c\xi-\frac{\xi^{2}}{\sigma}}\xi^{\frac{n+\beta-1}{2}}\left[\frac{n+\beta+1}{2}+\xi\left(c-\frac{2}{\sigma}\xi\right)\right]\\
        & = & 0
\end{eqnarray*}
iff
$$\frac{n+\beta+1}{2}+c\xi-\frac{2}{\sigma}\xi^{2}=0$$
iff
$$\xi=\frac{c\sigma+\sqrt{c^{2}\sigma^{2}+4\sigma(n+\beta+1)}}{4}=:\xi^{*}.$$
 So, $$\|f\|_{h}\leq 2^{-n-\frac{1+\beta}{4}}\pi^{-n-\frac{1}{4}}c^{\frac{1-n-\beta}{4}}\left[(\xi^{*})^{\frac{n+\beta+1}{2}}e^{c\xi^{*}-\frac{(\xi^{*})^{2}}{\sigma}}\right]^{1/2}\|f\|_{E_{\sigma}}.$$ We sum it up in the following theorem.
\begin{thm}
Under the conditions of Lemma 2.1, any $f\in E_{\sigma}$ satisfies
$$\|f\|_{h}\leq 2^{-n-\frac{1+\beta}{4}}\pi^{-n-\frac{1}{4}}c^{\frac{1-n-\beta}{4}}\left\{(\xi^{*})^{(n+\beta+1)/2}e^{c\xi^{*}-(\xi^{*})^{2}/\sigma}\right\}^{1/2}\|f\|_{E_{\sigma}}$$
where $$\xi^{*}:=\frac{c\sigma+\sqrt{c^{2}\sigma^{2}+4\sigma(n+\beta+1)}}{4}.$$
\end{thm}
\begin{cor}
Let $\sigma>0$ and $\beta<0$. If $|n+\beta|\geq 1$ and $n+\beta+1\geq 0$, then (4) in Theorem 1.3 has the form $$|f(x)-s(x)|\leq 2^{-\frac{3}{4}n}\pi^{-\frac{3}{4}n}\sqrt{n\alpha_{n}}\sqrt{\Delta_{0}}c^{\frac{1+\beta-n}{4}}\left\{(\xi^{*})^{(n+\beta+1)/2}e^{c\xi^{*}-(\xi^{*})^{2}/\sigma}\right\}^{1/2}(\lambda)^{\frac{1}{\delta}}\|f\|_{E_{\sigma}}$$
whenever $f\in E_{\sigma}$, where $\xi^{*}:=\frac{c\sigma+\sqrt{c^{2}\sigma^{2}+4\sigma(n+\beta+1)}}{4}$.
\end{cor}
{\bf Remark}. Note that Corollary 2.3 covers the most useful cases $\beta=-1$ and $n\geq2$. However the case $\beta=-1$ and $n=1$ is excluded. For $\beta=-1$ and $n=1$ we need a different approach.
\begin{lem}
Let $\sigma>0,\ \beta=-1$ and $n=1$. Then $E_{\sigma}\subseteq {\cal C}_{h,m}(R^{n})$ and for any $f\in E_{\sigma}$,
$$\|f\|_{h}\leq \frac{1}{2^{n+\frac{1}{4}}\pi^{n}}\left\{ \frac{1}{ln2}+2\sqrt{3}\sup_{|\xi|>\frac{1}{c}}\sqrt{c|\xi|}e^{c|\xi|-|\xi|^{2}/\sigma}\right\}^{1/2}\|f\|_{E_{\sigma}}.$$
\end{lem}
{\bf Proof}. Let $f\in E_{\sigma}$. By \cite{MN2} and \cite{Lu1},
\begin{eqnarray*}
\|f\|_{h} & = & \left\{ \frac{1}{(2\pi)^{2n}}\int\frac{|\hat{f}(\xi)|^{2}}{\hat{h}(\xi)}d\xi\right\}^{1/2}\\
          & = & \frac{1}{(2\pi)^{n}}\left\{\int\frac{|\hat{f}(\xi)|^{2}}{\sqrt{2}{\cal K}_{0}(c|\xi|)}d\xi\right\}^{1/2}\ (Theorem\ 8.15\ of\ \cite{We})\\
          & = & \frac{1}{(2\pi)^{n}2^{\frac{1}{4}}}\left\{\int_{|\xi|\leq\frac{1}{c}}\frac{|\hat{f}(\xi)|^{2}}{{\cal K}_{0}(c|\xi|)}d\xi +\int_{|\xi|>\frac{1}{c}}\frac{|\hat{f}(\xi)|^{2}}{{\cal K}_{0}(c|\xi|)}d\xi\right\}^{1/2}.
\end{eqnarray*}
Here,
\begin{eqnarray*}
\int_{|\xi|\leq\frac{1}{c}}\frac{|\hat{f}(\xi)|^{2}}{{\cal K}_{0}(c|\xi|)}dc\xi & \sim & \int_{|\xi|\leq \frac{1}{c}}\frac{|\hat{f}(\xi)|^{2}}{-ln\frac{c|\xi|}{2}}d\xi\ ({\cal K}_{0}(z)\sim -\{(ln\frac{z}{2})+r\}I_{0}(z)\ as\ z\rightarrow 0\ where\ r\sim 0.577\\
    &         & by\ p.\ 255,\ p.\ 374\ and\ p.\ 379\ of\ \cite{Ab})\\
    & = & \int_{|\xi|\leq\frac{1}{c}}|\hat{f}(\xi)|^{2}e^{\frac{|\xi|^{2}}{\sigma}}\frac{1}{e^{\frac{|\xi|^{2}}{\sigma}}|ln\frac{c|\xi|}{2}|}d\xi\\
    & \leq & \sup_{|\xi|\leq\frac{1}{c}}\left\{\frac{1}{e^{\frac{|\xi|^{2}}{\sigma}}|ln\frac{c|\xi|}{2}|}\right\}\cdot \int_{|\xi|\leq\frac{1}{c}}|\hat{f}(\xi)|^{2}e^{\frac{|\xi|^{2}}{\sigma}}d\xi\\
    & \leq & \frac{1}{e^{\frac{|\xi^{*}|^{2}}{\sigma}}|ln\frac{c|\xi^{*}|}{2}|}\cdot \|f\|^{2}_{E_{\sigma}}\ where\ 0<|\xi^{*}|\leq\frac{1}{c}\ for\ some\ \xi^{*}\\
    & \leq & \frac{1}{|ln2|}\cdot \|f\|^{2}_{E_{\sigma}}.
\end{eqnarray*}
Also, since $\Gamma(\frac{1}{2})=\sqrt{\pi}$,
\begin{eqnarray*}
\int_{|\xi|>\frac{1}{c}}\frac{|\hat{f}(\xi)|^{2}}{{\cal K}_{0}(c|\xi|)}d\xi & \leq & 2\sqrt{3}\int_{|\xi|>\frac{1}{c}}|\hat{f}(\xi)|^{2}e^{\frac{|\xi|^{2}}{\sigma}}\frac{\sqrt{c|\xi|}e^{c|\xi|}}{e^{\frac{|\xi|^{2}}{\sigma}}}d\xi\\
& \leq & 2\sqrt{3}\sup_{|\xi|>\frac{1}{c}}\frac{\sqrt{c|\xi|}e^{c|\xi|}}{e^{\frac{|\xi|^{2}}{\sigma}}}\cdot \|f\|^{2}_{E_{\sigma}}.
\end{eqnarray*}
Our lemma thus follows immediately. \hspace{9.3cm} $\sharp$
\begin{thm}
Let $\sigma>0,\ \beta=-1$ and $n=1$. For any $f\in E_{\sigma}$,
$$\|f\|_{h}\leq2^{-(n+\frac{1}{4})}\pi^{-n}\left\{\frac{1}{ln2}+2\sqrt{3}M(c)\right\}^{1/2}\|f\|_{E_{\sigma}}$$
where $M(c):=e^{1-\frac{1}{\sigma c^{2}}}$ if $c\leq\frac{2}{\sqrt{3\sigma}}$ and $M(c):=g\left(\frac{c\sigma+\sqrt{c^{2}\sigma^{2}+4\sigma}}{4}\right)$ if $c>\frac{2}{\sqrt{3\sigma}}$, where $g(\xi):=\sqrt{c\xi}e^{c\xi-\frac{\xi^{2}}{\sigma}}$.
\end{thm}
{\bf Proof}. The maximum of $g(\xi)$ on $[\frac{1}{c}, \infty)$ obviously exists. In order to find its exact value, we first find $g'(\xi)$. Note that
\begin{eqnarray*}
g'(\xi) & = & \frac{e^{\frac{\xi^{2}}{\sigma}}[D_{\xi}\sqrt{c\xi}e^{c\xi}]-\sqrt{c\xi}e^{c\xi}D_{\xi}e^{\frac{\xi^{2}}{\sigma}}}{e^{\frac{2\xi^{2}}{\sigma}}}\\
        & = & \frac{e^{\frac{\xi^{2}}{\sigma}}[\frac{c}{2\sqrt{c\xi}}e^{c\xi}+\sqrt{c\xi}e^{c\xi}c]-\sqrt{c\xi}e^{c\xi}e^{\frac{\xi^{2}}{\sigma}}\frac{2\xi}{\sigma}}{e^{\frac{2\xi^{2}}{\sigma}}}\\
        & = & \frac{\frac{c}{2\sqrt{c\xi}}e^{c\xi}+c\sqrt{c\xi}e^{c\xi}-\frac{2}{\sigma}\xi\sqrt{c\xi}e^{c\xi}}{e^{\frac{\xi^{2}}{\sigma}}}\\
        & = & \frac{e^{c\xi}\left[\frac{\sqrt{c}}{2\sqrt{\xi}}+c\sqrt{c\xi}-\frac{2}{\sigma}\xi\sqrt{c\xi}\right]}{e^{\frac{\xi^{2}}{\sigma}}}\\
        & = & 0
\end{eqnarray*}
iff
$$\frac{\sqrt{c}}{2\sqrt{\xi}}+c\sqrt{c\xi}-\frac{2}{\sigma}\xi\sqrt{c\xi}=0$$
iff
$$\xi^{2}-\frac{c\sigma}{2}\xi-\frac{\sigma}{4}=0$$
iff
$$\xi=\frac{c\sigma+\sqrt{c^{2}\sigma^{2}+4\sigma}}{4}.$$Let $\xi^{*}:=\frac{c\sigma+\sqrt{c^{2}\sigma^{2}+4\sigma}}{4}$. Then $\xi^{*}\leq\frac{1}{c}$ iff $c\leq\frac{2}{\sqrt{3\sigma}}$. Also, \(\lim_{\xi\rightarrow 0^{+}}g(\xi)=\lim_{\xi\rightarrow \infty}g(\xi)=0\). This gives that\\
$$\sup_{\xi>\frac{1}{c}}g(\xi)=\left\{ \begin{array}{ll}
                                        g(\frac{1}{c})=e^{1-\frac{1}{\sigma c^{2}}}        & \mbox{if $c\leq\frac{2}{\sqrt{3\sigma}},$}\\
                                        g(\xi^{*})=g\left(\frac{c\sigma+\sqrt{c^{2}\sigma^{2}+4\sigma}}{4}\right) & \mbox{if $c>\frac{2}{\sqrt{3\sigma}}.$}
                                      \end{array}
\right. $$The theorem then follows from Lemma 2.4. \hspace{8.4cm}  $\sharp$
\begin{cor}
Let $\sigma>0,\ \beta=-1$ and $n=1$. For any $f\in E_{\sigma}$, formula (4) in Theorem 1.3 can be expressed as
$$|f(x)-s(x)|\leq 2^{-1}\sqrt{\Delta_{0}}(\lambda)^{\frac{1}{\delta}}\frac{1}{\sqrt{c}}\left\{ \frac{1}{ln2}+2\sqrt{3}M(c)\right\}^{1/2}\|f\|_{E_{\sigma}}.$$
\end{cor}
Now we begin the study of the case $\beta>0$.
\begin{lem}
Let $\sigma>0,\ \beta>0$ and $n\geq 1$. Then $E_{\sigma}\subseteq{\cal C}_{h,m}(R^{n})$ and for any $f\in E_{\sigma}$,
$$\|f\|_{h}\leq d_{0}c^{\frac{1-\beta-n}{4}}\left\{ \sup_{\xi\in R^{n}}\frac{|\xi|^{(1+\beta+n)/2}e^{c|\xi|}}{e^{|\xi|^{2}/\sigma}}\right\}^{1/2}\|f\|_{E_{\sigma}}$$
where $d_{0}$ is a constant depending on $n,\beta$ only.
\end{lem}
{\bf Proof}. By definition,
\begin{eqnarray*}
  \|f\|_{h} & = & \left\{ \sum_{|\alpha|=m}\frac{m!}{\alpha!}\| (D^{\alpha}f)^{\wedge}\|^{2}_{L^{2}(\rho)}\right\}^{1/2}\ (by\ Corollary\ 3.3\ of\ \cite{MN2}\ and\ \cite{Lu1})\\
            & = & \left\{\sum_{|\alpha|=m}\frac{m!}{\alpha!}\int|(D^{\alpha}f)^{\wedge}(\xi)|^{2}d\rho\right\}^{1/2}\\
            & = & \left\{\sum_{|\alpha|=m}\frac{m!}{\alpha!}\int|(D^{\alpha}f)^{\wedge}(\xi)|^{2}\cdot \frac{1}{(2\pi)^{2n}|\xi|^{2m}\hat{h}(\xi)}d\xi\right\}^{1/2}\ (by\  \cite{MN2})\\
            & = & \left\{\sum_{|\alpha|=m}\frac{m!}{\alpha!}\int|i^{m}\xi^{\alpha}\hat{f}(\xi)|^{2}\cdot \frac{1}{(2\pi)^{2n}|\xi|^{2m}\hat{h}(\xi)}d\xi \right\}^{1/2}\\
            & = & \frac{(m!)^{\frac{1}{2}}}{(2\pi)^{n}}\left\{\sum_{|\alpha|=m}\frac{1}{\alpha!}\int\frac{\xi^{2\alpha}|\hat{f}(\xi)|^{2}}{|\xi|^{2m}\hat{h}(\xi)}d\xi\right\}^{1/2}\\
            & \leq & \frac{\sqrt{m!}}{(2\pi)^{n}\sqrt{2^{1+\frac{\beta}{2}}}}\left\{C(m,n)\int|\hat{f}(\xi)|^{2}\frac{|\xi|^{\frac{\beta+n}{2}}}{c^{\frac{\beta+n}{2}}{\cal K}_{\frac{n+\beta}{2}}(c|\xi|)}d\xi\right\}^{1/2}\ (by\ \cite{We})\\
            &   & where\ C(m,n)\ denotes\ the\ number\ of\ terms\ in\ \sum\\
            & \leq & \frac{\sqrt{m!C(m,n)}}{(2\pi)^{n}\sqrt{2^{1+\frac{\beta}{2}}}c^{\frac{\beta+n}{4}}}\left\{\int|\hat{f}(\xi)|^{2}|\xi|^{\frac{\beta+n}{2}}\cdot \frac{1}{\sqrt{\frac{\pi}{2}}\cdot \frac{e^{-c|\xi||}}{\sqrt{c|\xi|}}}d\xi\right\}^{1/2}\ (by\ \cite{We})\\
            & = &  \frac{\sqrt{m!C(m,n)}}{(2\pi)^{n}\sqrt{2^{1+\frac{\beta}{2}}}}\cdot c^{\frac{1-(\beta+n)}{4}}\left(\sqrt{\frac{2}{\pi}}\right)^{\frac{1}{2}}\left\{\int|\hat{f}(\xi)|^{2}|\xi|^{\frac{1+\beta+n}{2}}e^{c|\xi|}d\xi\right\}^{1/2}\\
            & \leq & d_{0}c^{\frac{1-\beta-n}{4}}\left\{\sup_{\xi\in R^{n}}\frac{|\xi|^{\frac{1+\beta+n}{2}}e^{c|\xi|}}{e^{\frac{|\xi|^{2}}{\sigma}}}\right\}^{1/2}\|f\|_{E_{\sigma}}\   where\ d_{0}:=\frac{\sqrt{m!C(m,n)}}{(2\pi)^{n}\sqrt{2^{1+\frac{\beta}{2}}}}\cdot \left(\frac{2}{\pi}\right)^{1/2}.       \\            
\end{eqnarray*}Since $\|f\|_{h}<\infty,\ f\in {\cal C}_{h,m}$. \hspace{10.5cm} $\sharp$
\begin{thm}
Under the conditions of Lemma 2.7,
$$\|f\|_{h}\leq d_{0}c^{\frac{1-\beta-n}{4}}\left\{\frac{(\xi^{*})^{(1+\beta+n)/2}e^{c\xi^{*}}}{e^{(\xi^{*})^{2}/\sigma}}\right\}^{1/2}\|f\|_{E_{\sigma}}$$
where $\xi^{*}=\frac{c\sigma+\sqrt{c^{2}\sigma^{2}+4\sigma(1+\beta+n)}}{4}$.
\end{thm}
{\bf Proof}. Let $g(x):=\frac{x^{\frac{1+\beta+n}{2}}e^{cx}}{e^{\frac{x^{2}}{\sigma}}},\ x>0$. Then
\begin{eqnarray*}
g'(x) & = & e^{-\frac{2x^{2}}{\sigma}}\left\{e^{\frac{x^{2}}{\sigma}}\left[\frac{1+\beta+n}{2}x^{\frac{\beta+n-1}{2}}e^{c x}+x^{\frac{1+\beta+n}{2}}e^{cx}c\right]-x^{\frac{1+\beta+n}{2}}e^{cx}e^{\frac{x^{2}}{\sigma}}\frac{2}{\sigma}x\right\}\\
      & = & e^{-\frac{x^{2}}{\sigma}}e^{c x}\left[\frac{1+\beta+n}{2}x^{\frac{\beta+n-1}{2}}+c x^{\frac{1+\beta+n}{2}}-x^{\frac{1+\beta+n}{2}}\frac{2}{\sigma}x\right]\\
      & = & e^{cx-\frac{x^{2}}{\sigma}}x^{\frac{\beta+n}{2}}\left[\frac{1+\beta+n}{2}x^{-\frac{1}{2}}+cx^{\frac{1}{2}}-x^{\frac{1}{2}}\frac{2}{\sigma}x\right]\\
      & = & 0
\end{eqnarray*}
iff
$$\frac{1+\beta+n}{2}\cdot \frac{1}{\sqrt{x}}+c\sqrt{x}-\frac{2}{\sigma}x\sqrt{x}=0$$
iff
$$\frac{1+\beta+n+2cx-\frac{4}{\sigma}x^{2}}{2\sqrt{x}}=0$$
iff
$$4x^{2}-2c\sigma x-\sigma(1+\beta+n)=0$$
iff
$$x=\frac{2c\sigma+\sqrt{4c^{2}\sigma^{2}+16\sigma(1+\beta+n)}}{8}$$
iff
$$x=\frac{c\sigma+\sqrt{c^{2}\sigma^{2}+4\sigma(1+\beta+n)}}{4}.$$Then the theorem follows immediately from the preceding lemma. \hspace{4.5cm} $\sharp$
\begin{cor}
Let $\sigma>0,\ \beta>0$ and $n\geq 1$. For any $f\in E_{\sigma}$, (4) in Theorem 1.3 can be expressed as
$$|f(x)-s(x)|\leq 2^{\frac{n+\beta+1}{4}}\pi^{\frac{n+1}{4}}\sqrt{n\alpha_{n}}\sqrt{\Delta_{0}}(\lambda)^{\frac{1}{\delta}}d_{0}c^{\frac{1+\beta-n}{4}}\left\{\frac{(\xi^{*})^{(1+\beta+n)/2}e^{c\xi^{*}}}{e^{(\xi^{*})^{2}/\sigma}}\right\}^{1/2}\|f\|_{E_{\sigma}}$$
where $d_{0}$ is defined as in Lemma 2.7, and $\xi^{*}$ is defined as in Theorem 2.8.
\end{cor}

\section{{\bf How to choose $c$ ?}}

\hspace{5mm}Theoretically $0<\lambda<1$ and $\delta$ can be arbitrarily small. Therefore $\lambda^{\frac{1}{\delta}}$ of (4) is very influential. The value of $\lambda$ had been unknown for a long time. Fortunately it's clarified in \cite{Lu3}. This is a breakthrough and makes it possible to assess the influence of $c$ on the error bound. However the value of $\lambda^{\frac{1}{\delta}}$ highly depends on whether $b_{0}$ in Theorem 1.3 is fixed or not. Hence we discuss it separately.
\subsection{{\bf $b_{0}$ fixed}}

\hspace{5mm}Let $b_{0}$ in Theorem 1.3 be fixed. Then the requirement $\delta \leq \delta_{0}$ forces $c\geq c_{0}$ where $c_{0}:=12\rho \sqrt{n} e^{2n\gamma_{n}}\gamma_{n}(m+1)\delta$ if we require $c\leq c_{1}$ defined below. Now,
$$C=\left\{ \begin{array}{ll}
              2\rho'\sqrt{n}e^{2n\gamma_{n}}  & \mbox{if $c\in [c_{0}, c_{1}],$}\\
              \frac{2}{3b_{0}}                & \mbox{if $c\in [c_{1}, \infty),$}      
            \end{array} \right. $$
where $c_{1}:=3b_{0}\rho\sqrt{n}e^{2n\gamma_{n}}$. Since $\lambda=\left(\frac{2}{3}\right)^{\frac{1}{6C\gamma_{n}}}$, we have
$$\lambda^{\frac{1}{\delta}}=\left\{\begin{array}{ll}
                                      e^{\eta(\delta)c}   & \mbox{if $c\in [c_{0},c_{1}],$}\\
                                      \left(\frac{2}{3}\right)^{\frac{b_{0}}{4\gamma_{n}\delta}}                                    & \mbox{if $c\in [c_{1},\infty),$}
                                    \end{array} \right. $$
where $\eta(\delta):=\frac{ln\frac{2}{3}}{12\rho\sqrt{n}e^{2n\gamma_{n}}\gamma_{n}\delta}$.

 It's easily seen that $\lambda^{\frac{1}{\delta}}$ is a continuous function of $c$ and is independent of c whenever $c\geq c_{1}$.

Obviously $\lambda^{\frac{1}{\delta}}$ is influential only when $\delta$ is very small. The number $\gamma_{n}$ grows very fast as $n$ increases. Therefore $\lambda^{\frac{1}{\delta}}\approx 1$ for high dimensions; unless $\delta$ is extremely small.

With these understandings we can now begin our theoretical analysis of the optimal $c$.\\
\\
{\bf Case 1}. \fbox{$\beta=-1$ and $n=1$} Let $f\in E_{\sigma}$ and $h$ be the map defined in (1) with $\beta=-1$ and $n=1$. Corollary 2.6 shows that $$|f(x)-s(x)|\leq MN(c)\cdot \|f\|_{E_{\sigma}}$$ where
$$MN(c):=\left\{ \begin{array}{ll}
                         \frac{1}{2}\cdot (\frac{2}{3})^{\frac{c}{24\delta e^{4}}}\cdot \frac{1}{\sqrt{c}}\cdot \{ \frac{1}{ln2}+2\sqrt{3}M(c)\}^{1/2} & \mbox{if $24\delta e^{4}\leq c\leq 3b_{0}e^{4},$}\\
                        \frac{1}{2}\cdot (\frac{2}{3})^{\frac{b_{0}}{8\delta}}\cdot \frac{1}{\sqrt{c}}\cdot \{ \frac{1}{ln2}+ 2\sqrt{3}M(c)\} & \mbox{if $c\geq 3b_{0}e^{4},$}
                       \end{array} \right.  $$
where $b_{0}$ is the side length of the domain cube and $\delta$ is the fill distance defined in Theorem 1.3, and 
$$M(c):=\left\{ \begin{array}{ll}
                      e^{1-\frac{1}{\sigma c^{2}}}  & \mbox{if $c\leq \frac{2}{\sqrt{3\sigma}},$}\\
                     g(\frac{c\sigma +\sqrt{c^{2}\sigma^{2}+4\sigma}}{4})  & \mbox{if $c>\frac{2}{\sqrt{3\sigma}},$}  
                        \end{array} \right.   $$ 
where $g(\xi):=\sqrt{c\xi}e^{c\xi-\frac{\xi^{2}}{\sigma}}$. The value $c$ minimizing $MN(c)$ in the interval $[24\delta e^{4},\infty)$ is then the optimal choice of $c$.\\

The restriction $c\geq 24\delta e^{4}$,  which is just $c_{0}$ defined in the beginning of this subsection, is a drawback. However, empirical results show that the optimal value of $c$ never lies in the interval $(0, 24\delta e^{4})$. It's the same for Cases 2 and 3. Hence we essentially have dealt with all positive c in the interval $(0,\infty)$. We call $MN(c)$ the MN function and its graph the MN curve, in honor of Professors W. R. Madych and S. A. Nelson for their outstanding contributions to the development of RBF. Some examples are given in Figs. 1-5, where $N_{d}$ denotes the number of data points used.\\
\\
\begin{figure}[t]
\centering
\includegraphics[scale=1.0]{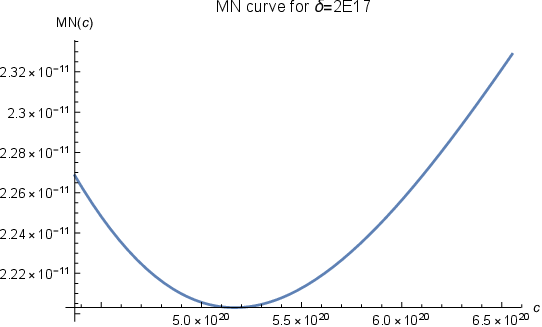}
\caption{Here $n=1,\beta=-1,\sigma=1$E$-42, b_{0}=1$E$18$ and $N_{d}=4$.}

\includegraphics[scale=1.0]{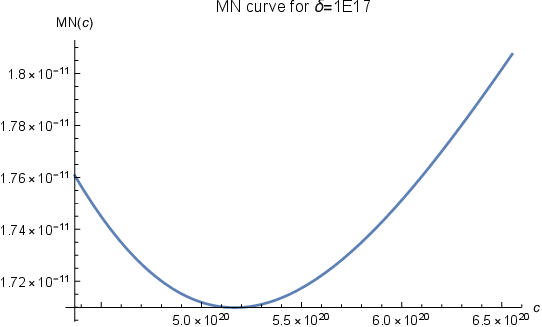}
\caption{Here $n=1,\beta=-1,\sigma=1$E$-42, b_{0}=1$E$18$ and $N_{d}=6$.}

\includegraphics[scale=1.0]{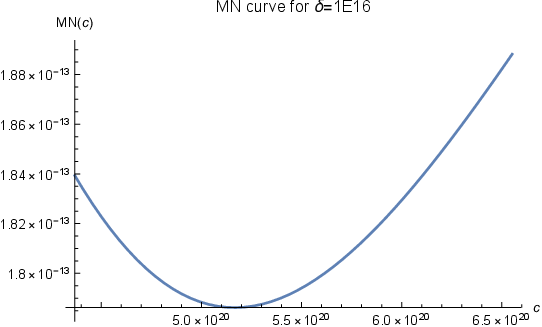}
\caption{Here $n=1,\beta=-1,\sigma=1$E$-42, b_{0}=1$E$18$ and $N_{d}=51$.}

\end{figure}

\clearpage

\begin{figure}[t]
\centering
\includegraphics[scale=1.0]{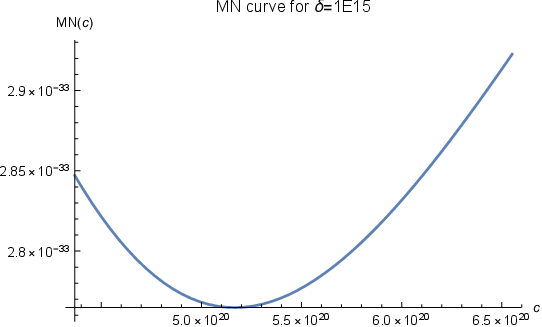}
\caption{Here $n=1,\beta=-1,\sigma=1$E$-42, b_{0}=1$E$18$ and $N_{d}=501$.}

\includegraphics[scale=1.0]{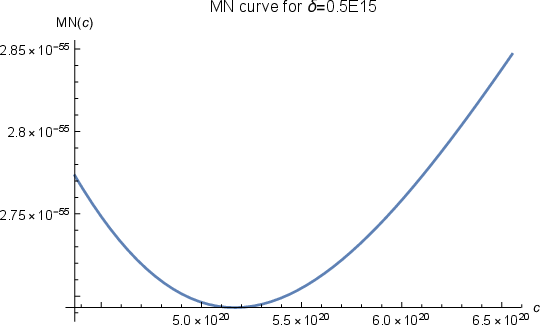}
\caption{Here $n=1,\beta=-1,\sigma=1$E$-42, b_{0}=1$E$18$ and $N_{d}=1001$.}

\end{figure}

\hspace{-5mm}{\bf Case 2}. \fbox{$|n+\beta|\geq 1$, $\beta<0$ and $n+\beta+1\geq 0$} Let $f\in E_{\sigma}$ and $E$ be the domain cube in Theorem 1.3 with side length $b_{0}$. Let $h$ be the map defined in (1) with $1+\beta+n\geq 0$, $\beta<0$ and $|n+\beta|\geq 1$. Corollary 2.3 shows that 
$$|f(x)-s(x)|\leq MN(c)\cdot \|f\|_{E_{\sigma}}$$
where  
$$MN(c):=\left\{ \begin{array}{ll}
                         (2\pi)^{-\frac{3}{4}n}\sqrt{n\alpha_{n}\Delta_{0}}c^{\frac{1+\beta-n}{4}}\{(\xi^{*})^{\frac{n+\beta+1}{2}}e^{c\xi^{*}-\frac{(\xi^{*})^{2}}{\sigma}}\}^{1/2}(\frac{2}{3})^{\frac{c}{12\rho \delta \gamma_{n}\sqrt{n}e^{2n\gamma_{n}}}} &   \mbox{if $c\in [c_{0}, c_{1}],$}\\
                         (2\pi)^{-\frac{3}{4}n}\sqrt{n\alpha_{n}\Delta_{0}}c^{\frac{1+\beta-n}{4}}\{(\xi^{*})^{\frac{n+\beta+1}{2}}e^{c\xi^{*}-\frac{(\xi^{*})^{2}}{\sigma}}\}^{1/2}(\frac{2}{3})^{\frac{b_{0}}{4\gamma_{n}\delta}} & \mbox{if $c\in [c_{1},\infty),$}
                         \end{array} \right.   $$ 
where $\delta$ is the fill distance, $\alpha_{n}$ is the volume of the unit ball in $R^{n}$, and $\gamma_{n}$, and $\Delta_{0}$, $\rho$ were defined in Defintions 1.1 and 1.2, respectively. The optimal value of $c$ is just the value minimizing $MN(c)$ in $[c_{0},\infty)$.
The MN curves of Case 2 are given in Figs. 6-10.\\
\\
\begin{figure}[t]
\centering
\includegraphics[scale=1.0]{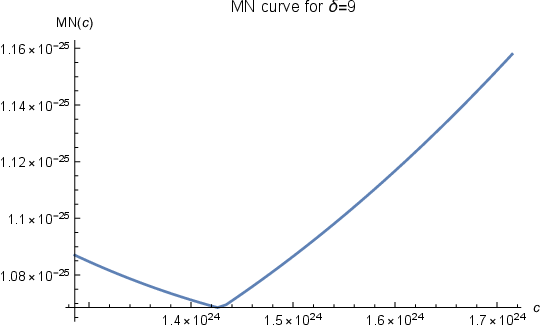}
\caption{Here $n=2,\beta=-1,\sigma=1$E$-48, b_{0}=480$ and $N_{d}=784$.}

\includegraphics[scale=1.0]{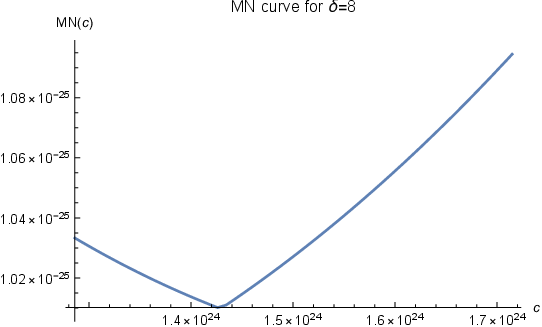}
\caption{Here $n=2,\beta=-1,\sigma=1$E$-48, b_{0}=480$ and $N_{d}=961$.}

\includegraphics[scale=1.0]{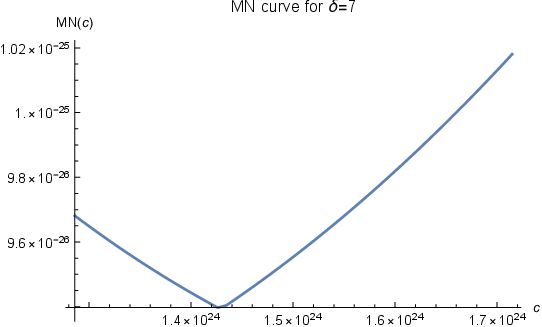}
\caption{Here $n=2,\beta=-1,\sigma=1$E$-48, b_{0}=480$ and $N_{d}=1296$.}

\end{figure}

\clearpage

\begin{figure}[t]
\centering
\includegraphics[scale=1.0]{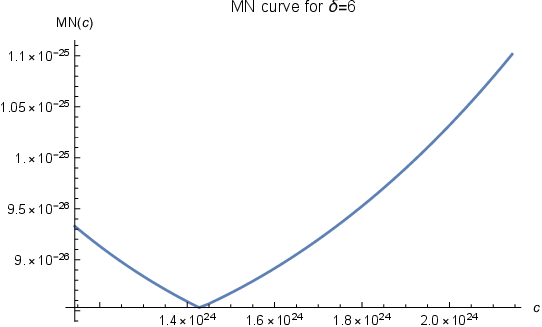}
\caption{Here $n=2,\beta=-1,\sigma=1$E$-48, b_{0}=480$ and $N_{d}=1681$.}

\includegraphics[scale=1.0]{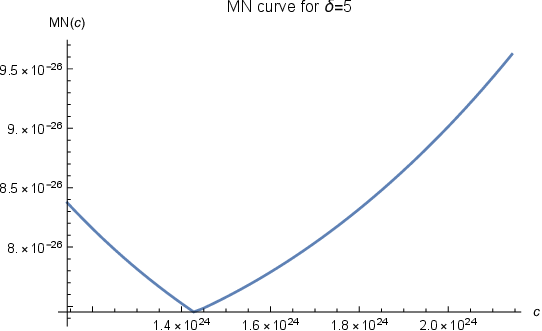}
\caption{Here $n=2,\beta=-1,\sigma=1$E$-48, b_{0}=480$ and $N_{d}=2041$.}

\end{figure}

\hspace{-5mm}{\bf Case 3}. \fbox{ $\beta>0$ and $n\geq 1$} Let $f\in E_{\sigma}$ and $E$ be the domain cube in Theorem 1.3 with side length $b_{0}$. Let $h$ be the map defined in (1) with $\beta>0$ and $n\geq 1$. Corollary 2.9 shows that 
$$|f(x)-s(x)|\leq d_{0}\cdot MN(c)\cdot \|f\|_{E_{\sigma}}$$
where
$$MN(c):=\left\{ \begin{array}{ll}
                         2^{\frac{n+\beta+1}{4}}\pi^{\frac{n+1}{4}}\sqrt{n\alpha_{n}\Delta_{0}}c^{\frac{1+\beta-n}{4}}\{\frac{(\xi^{*})^{(1+\beta+n)/2}e^{c\xi^{*}}}{e^{(\xi^{*})^{2}/\sigma}}\}^{1/2}(\frac{2}{3})^{\frac{c}{12\rho \sqrt{n}\gamma_{n}\delta e^{2n\gamma_{n}}}}  & \mbox{if  $c\in [c_{0}, c_{1}],$}\\
                         2^{\frac{n+\beta+1}{4}}\pi^{\frac{n+1}{4}}\sqrt{n\alpha_{n}\Delta_{0}}c^{\frac{1+\beta-n}{4}}\{\frac{(\xi^{*})^{(1+\beta+n)/2}e^{c\xi^{*}}}{e^{(\xi^{*})^{2}/\sigma}}\}^{1/2}(\frac{2}{3})^{\frac{b_{0}}{4\gamma_{n}\delta}}     & \mbox{if $c\in [c_{1},\infty),$}
                         \end{array} \right.   $$ 
where  $\delta$ is the fill distance, $\alpha_{n}$ is the volume of the unit ball in $R^{n}$, and $\gamma_{n}$, and $\Delta_{0}$, $\rho$ were defined in Defintions 1.1 and 1.2, respectively. The constant $d_{0}$ was defined in Lemma 2.7 and depends on $n$ and $\beta$ only. Then we choose $c$ which minimizes $MN(c)$ in $[c_{0},\infty)$. The MN curves appear in Figs. 11-15.

\begin{figure}[t]
\centering
\includegraphics[scale=1.0]{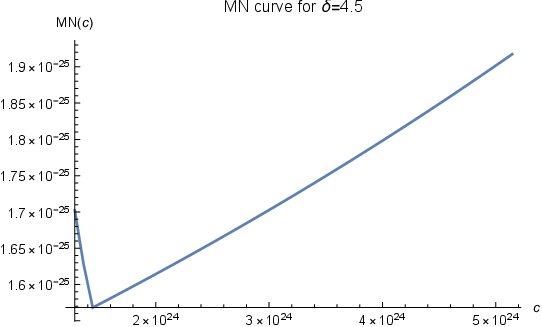}
\caption{Here $n=2,\beta=1,\sigma=1$E$-50, b_{0}=480$ and $N_{d}=3025$.}

\includegraphics[scale=1.0]{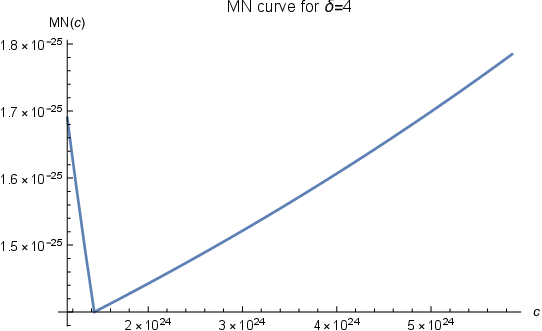}
\caption{Here $n=2,\beta=1,\sigma=1$E$-50, b_{0}=480$ and $N_{d}=3721$.}

\includegraphics[scale=1.0]{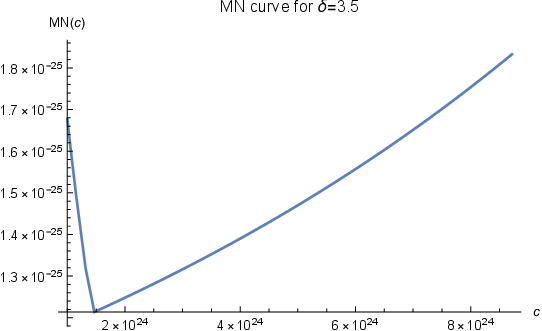}
\caption{Here $n=2,\beta=1,\sigma=1$E$-50, b_{0}=480$ and $N_{d}=4900$.}

\end{figure}

\clearpage

\begin{figure}[t]
\centering
\includegraphics[scale=1.0]{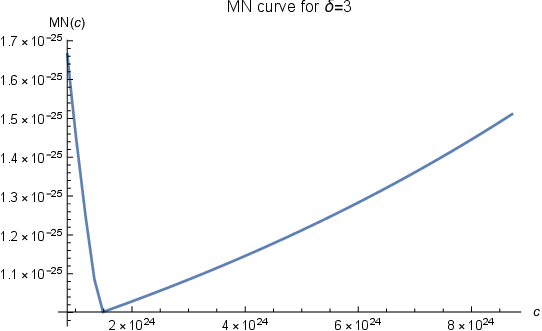}
\caption{Here $n=2,\beta=1,\sigma=1$E$-50, b_{0}=480$ and $N_{d}=6561$.}

\includegraphics[scale=1.0]{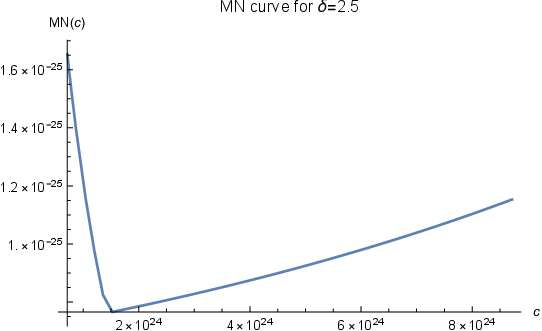}
\caption{Here $n=2,\beta=1,\sigma=1$E$-50, b_{0}=480$ and $N_{d}=9409$.}

\end{figure}

Something amazing in Cases 2 and 3 is that increasing the number of data points does not improve the MN function values much. What's important is the choice of $c$. Although in Case 1 the essential error bound $MN(c)$ for the interpolation greatly depends on the fill distance $\delta$, very small $MN(c)$ can be achieved by using very few data points. For example, in Fig 1, only 4 data points are used in a huge domain of length 1E18. In fact, this kind of results is frequently seen in our $c$ theory, as long as $c$ is chosen well.

\subsection{{\bf $b_{0}$ not fixed}}

\hspace{5mm}As explained in \cite{Lu4} and \cite{MN3}, some domains are invariant under dilation. Any point in such a domain is contained in a cube of side $b_{0}$ where $b_{0}$ can be made arbitrarily large and the cube is still contained in the domain. For example,
$$\Omega:=\left\{(x_{1},\cdots,x_{n}):\ x_{i}\geq 0\ for\ i=1,\cdots,n\right\}$$
is such a domain. So is $\Omega=R^{n}$.

In Theorem 1.3, if $b_{0}$ can be made arbitrarily large, then both $C$ and $\lambda$ can be made arbitrarily small by increasing $c$ and $b_{0}$. The optimal choice of $c$ will hence be very different.

In this paper every approximated function $f$ belongs to $E_{\sigma}$. The domain of $f$ is of course $R^{n}$. However the interpolation occurs in a cube as required in Theorem 1.3. In practical problems the interpolation often can occur only in some subset $\Omega$ of $R^{n}$, even if the domain of the approximated function is the entire $R^{n}$. In this subsection the dilation-invariant domain $\Omega$ denotes the subset of $R^{n}$ where interpolation can occur.

We begin with the case $\beta=-1$ and $n=1$.\\
\\
{\bf Case 1}. \fbox{$\beta=-1$ and $n=1$} Let $f\in E_{\sigma}$ be the interpolated function and $\Omega\subseteq R^{1}$ be such that for any $x\in\Omega$ and $b_{0}>0$, there exists a cube $E$ of side $b_{0}$ such that $x\in E\subseteq \Omega$ and interpolation can occur in $E$. Let $h$ be the map defined in (1) with $\beta=-1$ and $n=1$. For any $\delta>0$, the optimal choice of $c$ in the interval $[24\delta e^{4},\infty)$ is the number minimizing $MN(c)$ where $MN(c)$ was defined in Case 1 of section 3 for $c\in [24\delta e^{4}, 3b_{0}e^{4}]$ .\\
\\
{\bf Reason}:  Since $C$ in Theorem 1.3 can be kept equal to $2\rho'\sqrt{n}e^{2n\gamma_{n}}$ by increasing $b_{0}$. \hspace{10.1cm} $\sharp$\\

What's noteworthy is that we only increase $b_{0}$ to keep $C=2\rho'\sqrt{n}e^{2n\gamma_{n}}$. We never decrease $b_{0}$ because it will only increase $C$ and make $\lambda$ and $\delta_{0}$ in Theorem 1.3 worse. With the same reason as Case 1, we have the following two cases.\\
\\
{\bf Case 2}. \fbox{$|n+\beta|\geq 1$, $\beta<0$ and $n+\beta+1\geq 0$} Let $f\in E_{\sigma}$ and be the interpolated function and $\Omega\subseteq R^{n}$ be such that for any $x\in\Omega$ and $b_{0}>0$, there exists a cube $E$ of side $b_{0}$ such that $x\in E\subseteq \Omega$ and interpolation can occur in $E$. Let $h$ be the map defined in (1) with $1+\beta+n\geq 0$, $\beta<0$ and $|n+\beta|\geq 1$. Then for any fixed $\delta>0$, the optimal choice of $c$ in the interval $[c_{0},\infty)$ is the number minimizing $MN(c)$ where $MN(c)$ was defined in Case 2 of Section 3. for $c\in [c_{0},c_{1}]$.\\
\\
{\bf Case 3}. \fbox{ $\beta>0$ and $n\geq 1$} Let $f\in E_{\sigma}$ and be the interpolated function and $\Omega\subseteq R^{n}$ be such that for any $x\in\Omega$ and $b_{0}>0$, there exists a cube $E$ of side $b_{0}$ such that $x\in E\subseteq \Omega$ and interpolation can occur in $E$. Let $h$ be the map defined in (1) with $\beta > 0$ and $n\geq 1$. Then for any fixed $\delta>0$, the optimal choice of $c$ in the interval $[c_{0},\infty)$ is the number minimizing $MN(c)$ where $MN(c)$ was defined in Case 3 of Section 3. for $c\in [c_{0},c_{1}]$.

\end{document}